\documentclass{article}
\usepackage{amssymb}
\usepackage{graphicx}
\usepackage{amsmath}

\newtheorem{theorem}{Theorem}

\newtheorem{conjecture}[theorem]{Conjecture}
\newtheorem{corollary}[theorem]{Corollary}

\newtheorem{lemma}[theorem]{Lemma}

\newtheorem{proposition}[theorem]{Proposition}

\begin{document}

\title{Rank gradient, cost of groups and the rank versus Heegaard genus problem}
\author{Mikl\'{o}s Ab\'{e}rt and Nikolay Nikolov}
\maketitle

\begin{abstract}
We study the growth of the rank of subgroups of finite index in residually
finite groups, by relating it to the notion of cost.

As a by-product, we show that the `Rank vs. Heegaard genus' conjecture on
hyperbolic $3$-manifolds is incompatible with the `Fixed Price problem' in
topological dynamics.
\end{abstract}

\section{Introduction}

Let $\Gamma $ be a finitely generated group. A \emph{chain} in $\Gamma $ is
a sequence $\Gamma =\Gamma _{0}\geq \Gamma _{1}\geq \ldots $ of subgroups of
finite index in $\Gamma $. Let $T=T(\Gamma ,(\Gamma _{n}))$ denote the \emph{%
coset tree} of the chain, a rooted tree on the set of right cosets of the
subgroups $\Gamma _{n}$ with edges $(\Gamma _{n}g,\Gamma _{n+1}g)$ for all $%
g\in \Gamma $ and $n\in \mathbb{N}$. The \emph{boundary} $\partial T$ of $T$
is the set of infinite rays starting from the root; it is naturally endowed
with the product topology and product measure coming from the tree. The
group $\Gamma $ acts by automorphisms on $T$; this action extends to measure
preserving homeomorphisms of the boundary.

We say that a chain $(\Gamma _{n})$ is \emph{Farber}, if the action of $%
\Gamma $ on the boundary of its coset tree $T=T(\Gamma ,(\Gamma _{n}))$ is
essentially free, that is, if almost every element of $\partial T$ has
trivial stabilizer in $\Gamma $. This is the case for example when the chain
consists of normal subgroups of $\Gamma $ and their intersection is trivial.
Note that then $\partial T$ is simply the profinite completion of $\Gamma $
with respect to $(\Gamma _{n})$ endowed with the normalized Haar measure.

For a group $G$ let $d(G)$ denote the minimal number of generators (or rank)
of $G$. Let the \emph{rank gradient} of $\Gamma $ with respect to $(\Gamma
_{n})$ be defined as 
\begin{equation*}
\mathrm{RG}(\Gamma ,(\Gamma _{n}))=\lim_{n\rightarrow \infty }\frac{d(\Gamma
_{n})-1}{\left| \Gamma :\Gamma _{n}\right| }
\end{equation*}
This notion has been introduced by Lackenby \cite{lack}.

Our first theorem relates the rank gradient of a Farber chain to the cost of
the action of the group on the boundary of the coset tree. The analytic
notion of cost was introduced by Levitt \cite{levit} and used by Gaboriau 
\cite{gabor} to show that free groups of different rank do not admit orbit
equivalent measurable actions (see also the book of Kechris and Miller \cite
{kechmill}).

\begin{theorem}
\label{coststable}Let $(\Gamma _{n})$ be a Farber chain in $\Gamma $. Then 
\begin{equation*}
\mathrm{RG}(\Gamma ,(\Gamma _{n}))=\mathrm{cost}(E)-1
\end{equation*}
where $E$ denotes the orbit relation given by the action of $\Gamma$ on the
boundary of the coset tree $T(\Gamma ,(\Gamma _{n}))$.
\end{theorem}

Theorem \ref{coststable} allows us to clash two well-known problems, one in $%
3$-manifold theory and the other in topological dynamics.

\medskip

\noindent \textbf{Rank vs Heegaard genus conjecture.} \emph{Let }$M$\emph{\
be a compact, orientable, hyperbolic }$3$\emph{-manifold. Then the Heegaard
genus of }$M$\emph{\ equals the rank of the fundamental group of }$M$\emph{.}

\medskip

It is easy to see that the Heegaard genus is always greater or equal to the
rank. The problem dates back to Waldhausen \cite{W}, who asked it for
arbitrary $3$-manifolds. This was proved false for Seifert manifolds by
Boileau and Zieschang in \cite{BZ} (see also \cite{schweid}), but it
remained open for hyperbolic $3$-manifolds. For the above formulation, see 
\cite[Conjecture 1.1]{shalen}. Also, it is not known, whether the ratio of
the two quantities can become arbitrarily large even for arbitrary $3$%
-manifolds; the best known lower bound comes from the Boileau-Zieschang
result.

\medskip

A countable group $\Gamma $ has \emph{fixed price}, if every essentially
free measure-preserving Borel action of $\Gamma $ has the same cost.
Gaboriau \cite{gabor} established fixed price for a large class of groups,
including free groups, higher rank non-uniform real irreducible lattices and
groups containing an infinite amenable normal subgroup and asked whether the
following holds.

\medskip

\noindent \textbf{Fixed Price problem.} \emph{Does every countable group
have fixed price?} \medskip

\begin{theorem}
\label{contradiction}Either the Rank vs Heegaard genus conjecture is false
or the Fixed price problem has a negative solution.
\end{theorem}

Moreover, if the Fixed Price problem has an affirmative answer, then the
Rank vs Heegaard genus conjecture fails in the following strong senses.
First, the ratio of the Heegaard genus and the rank of the fundamental group
of a compact, orientable hyperbolic $3$-manifold can get arbitrarily large.
Second, the counterexamples are not exotic, rather this seems to be the
general asymptotic behaviour of arithmetic hyperbolic $3$-manifolds.

\medskip

The contradiction between the two problems is established via the (unknown)
answer to the following question: Does the rank gradient $\mathrm{RG}%
(\Gamma, (\Gamma_n))$ depend on the choice of the Farber chain $(\Gamma_n)$
in $\Gamma$?

If it does, then Theorem \ref{coststable} trivially provides a negative
answer to the Fixed Price problem. In the other direction, there are
specific (uniform or non-uniform) arithmetic lattices in $\mathrm{SL}_{2}(%
\mathbb{C})$ (e.g. $\mathrm{SL}_{2}(\mathbb{Z}[i])$), that viewed as
abstract groups, possess a chain of subgroups with vanishing rank gradient.
On the other hand, being arithmetic groups, they have property ($\tau $)
with respect to congruence subgroups. Using work of Lackenby \cite{lack},
this allows one to construct a covering tower of the $3$-manifold
corresponding to the lattice, where the Heegaard genus grows linearly. Now,
if the rank gradient is independent of the chain, the rank must grow
sublinearly on this tower, making the ratio of the Heegaard genus and the
rank arbitrarily large.

\medskip

We shall derive the independence of the rank gradient from the chain from a
hypothesis that is much weaker than fixed price (see \cite[Problem 25.5]
{kechmill}).

\medskip

\noindent \textbf{Multiplicativity of cost-1 problem.} \emph{Let }$\Gamma $%
\emph{\ be a measurable, essentially free action on }$(X,\mu )$\emph{\ and
let }$H$\emph{\ be a subgroup of }$\Gamma $\emph{\ of finite index. Is it
true that } 
\begin{equation*}
\mathrm{cost}(H,X)-1=(\mathrm{cost}(\Gamma ,X)-1)\left| \Gamma :H\right| 
\text{?}
\end{equation*}

That is, does the cost of an action behave the same way as the rank for free
groups?

\medskip

The whole theory developed in this paper has a close connection to $L^{2}$
Betti numbers. The L\"{u}ck approximation result \cite{luck} implies that
for finitely presented groups and normal chains with trivial intersection,
if we replace $d(\Gamma _{n})$ with the first Betti number $\beta
_{1}(\Gamma _{n})$ in the definition of rank gradient, then the limit equals
the first $L^{2}$ Betti number $\beta _{1}^{2}(\Gamma )$. This has been
generalized by Farber \cite{farber} to chains satisfying his condition (see
also \cite{berggab} for examples showing the necessity of this condition).
On the other hand, Gaboriau \cite{gabihes} has introduced $L^{2}$ Betti
numbers of a measurable equivalence relation $E$ and asked whether $\beta
_{1}^{2}(E)=\mathrm{cost}(E)-1$ in general. An affirmative answer to
Gaboriau's question would imply the surprising result that the asymptotic
growth of $\beta _{1}(\Gamma _{n})$ and $d(\Gamma _{n})$ are equal for
arbitrary groups.

\medskip

Theorem \ref{coststable} immediately allows us to compute the rank gradient
for a class of groups where elementary methods seem to fail working.

\begin{theorem}
\label{amennormal}Let $\Gamma $ be a residually finite group with an
infinite amenable normal subgroup and let $\Gamma _{n}$ be a Farber chain in 
$\Gamma$. Then $\mathrm{RG}(\Gamma ,(\Gamma _{n}))=0$.
\end{theorem}

This generalizes a result of Lackenby \cite{lack} who proved the result for
finitely presented amenable groups.

We also answer a question of Kechris and Miller \cite[Problem 35.7]{kechmill}%
. They asked whether for a countable infinite Euclidean domain $D$ the group 
$\mathrm{SL}(2,D)$ has cost $1$ if and only if $D$ has infinitely many
units. The answer is negative: as we shall see, $\mathrm{SL}(2,\mathbb{Z}%
[i]) $ has cost $1$ and $\mathbb{Z}[i]$ has finitely many units. A
well-known conjecture by Thurston asserts that every hyperbolic $3$-manifold
virtually fibers over the circle. This would imply that any lattice in $%
\mathrm{SL}(2,\mathbb{C})$ has cost $1$. On the other hand, as we show in
Section \ref{manifolds}, lattices in $\mathrm{SL}(2,\mathbb{R})$ have
positive rank gradient with respect to any Farber chain.

\medskip

The paper is organized as follows. In Section \ref{basics} we define
measurable actions, graphings and the cost. We also introduce a new type of
graphing and an invariant, the product cost. We relate the product cost to
the cost and prove a general approximation result on graphings. In Section 
\ref{prodcost} we apply this result to chains and express the product cost
in terms of the rank gradient, thus proving Theorem \ref{coststable}. In
Section \ref{examples} we use the theory developed to prove Corollary \ref
{amennormal} and compute the rank gradient for some important classes of
groups. We also discuss what happens if we relax the Farber condition on the
chain. Finally, in Section \ref{manifolds} we consider the rank gradient of
lattices in Lie groups, prove Theorem \ref{contradiction} and answer the
question of Kechris and Miller.

\medskip

\noindent \textbf{Acknowledgements.} The authors thank Ian Agol, Gabor Elek, Marc
Lackenby and Alan Reid for helpful comments.

\section{Graphings and cost \label{basics}}

In this section we introduce Borel actions and the standard notion of cost
as well as a new version of cost that will be useful in proving Theorem \ref
{coststable}.

Let $X$ be a standard Borel space. Let the countable group $\Gamma $ act on $%
X$ by Borel-automorphisms. Let $\mu $ be a $\Gamma $-invariant probability
measure on $X$. We assume that the $\Gamma $-action has finitely many
ergodic components.

Let us define the relation $E$ on $X$ by 
\begin{equation*}
xEy\text{ if there exists }\gamma \in \Gamma \text{ with }y=x\cdot \gamma
\end{equation*}
Then $E$ is a Borel equivalence relation and every equivalence class is
countable.

Since $E$ is a subset of $X\times X$, it can also be considered as a graph
on $X$. A \emph{Borel subgraph} of $E$ is a directed graph on $X$ such that
the edge set is a Borel subset of $E$.

Let $S\subseteq X\times X$ be an arbitrary graph on $X$. A path from $x$ to $%
y$ in $S$ of length $k\geq 1$ is a sequence $x_{0},x_{1},\dots ,x_{k}\in X$
such that:

\begin{itemize}
\item  $x_{0}=x$, $x_{k}=y$;

\item  $(x_{i},x_{i+1})\in S$ or $(x_{i+1},x_{i})\in S$ ($0\leq i\leq k-1$).
\end{itemize}

Note that we consider undirected paths above. For $k\geq 1$ let us define
the graph $S^{k}$ by 
\begin{equation*}
(x,y)\in S^{k}\text{ if }x=y\text{ or there is a path from }x\text{ to }y%
\text{ in }S\text{ of length at most }k
\end{equation*}

We say that a subgraph $S$ of $E$ \emph{spans} $E$, if for any $(x,y)\in E$
with $x\neq y$ there exists a path from $x$ to $y$ in $S$. Trivially, this
holds if and only if 
\begin{equation*}
\bigcup_{n}S^{n}=E
\end{equation*}

The \emph{edge-measure} of a Borel subgraph $S$ of $E$ is defined as 
\begin{equation*}
e(S)=\int_{x\in X}\deg _{S}(x)\,d\mu
\end{equation*}
where $\deg (x)$ is the number of edges in $S$ with initial vertex $x$: 
\begin{equation*}
\deg _{S}x=\left| \left\{ y\in X\mid (x,y)\in S\right\} \right|
\end{equation*}
Note that $e(S)$ may be infinite. The \emph{cost} of $E$ is defined as 
\begin{equation*}
\mathrm{cost}(E)=\mathrm{cost}(\Gamma ,X)=\inf e(S)
\end{equation*}
where the infimum is taken over all Borel subgraphs $S$ of $E$ that span $E$%
. The cost of $\Gamma $ is defined as 
\begin{equation*}
\mathrm{cost}(\Gamma )=\inf \mathrm{cost}(\Gamma ,X)
\end{equation*}
where the infimum is taken over all ergodic, essentially free actions of $%
\Gamma $ on a standard Borel space $X$.

We say that $\Gamma $ has \emph{fixed price} $c$ if all ergodic, essentially
free actions of $\Gamma $ on a standard Borel space $X$ have cost $c$. It is
not known whether every countable group has fixed price \cite{gabor}.

If $\Gamma $ is generated by $g_{1},\ldots ,g_{d}$ then it is easy to see
that the set 
\begin{equation*}
\bigcup_{i=1}^{d}\bigcup_{x\in X}\left\{ (x,x\cdot g_{i})\right\}
\end{equation*}
is a spanning Borel subrelation of $E$ of edge measure $d$. This implies the
following.

\begin{lemma}
\label{kisranglemma}$\mathrm{cost}(\Gamma ,X)\leq d(\Gamma )$.
\end{lemma}

Now we will look at cost from another point of view, using graphings. The
notion basically comes from \cite{levit} but we will need to use it somewhat
differently. The advantage is that this notion will work for
non-essentially-free (even finite) actions as well.

Let us consider the product space $X\times \Gamma $ where $\Gamma $ is
endowed with the discrete topology and the counting measure. Denote the
product measure by $e$. A \emph{graphing} is a Borel subset of $X\times
\Gamma $. For a graphing $M$ and $\gamma \in \Gamma $ let 
\begin{equation*}
M_{\gamma }=\left\{ x\in X\mid (x,\gamma )\in M\right\}
\end{equation*}
be the $\gamma $-fiber of $M$. For $x\in X$ let 
\begin{equation*}
N(M,x)=\left\{ \gamma \in \Gamma \mid x\in M_{\gamma }\right\} \text{ and }%
\deg _{M}x=\left| N(M,x)\right|
\end{equation*}
be the set of neighbours and the degree of $x$ in $M$. Using this notation,
we have 
\begin{equation*}
e(M)=\sum_{\gamma \in \Gamma }\mu (M_{\gamma })=\int_{x\in X}\deg
_{M}(x)\,d\mu
\end{equation*}

We will need a definition of powering of graphings. Let $I$ be the graphing
defined by 
\begin{equation*}
I_{\gamma }=\left\{ 
\begin{array}{cc}
X & \gamma =1 \\ 
\emptyset & \text{otherwise}
\end{array}
\right.
\end{equation*}
For a graphing $M$ let the graphing $M^{\intercal }$ be defined by 
\begin{equation*}
M_{\gamma }^{\intercal }=M_{\gamma ^{-1}}\cdot \gamma ^{-1}\text{ (}\gamma
\in \Gamma \text{)}
\end{equation*}
and let $\overline{M}$ be defined by 
\begin{equation*}
\overline{M}=M\cup M^{\intercal }\cup I
\end{equation*}
For graphings $M$ and $N$ let us define the graphing $M\cdot N$ by 
\begin{equation*}
(M\cdot N)_{\gamma }=\bigcup_{\delta \in \Gamma }\left( M_{\delta }\cap
(N_{\delta ^{-1}\gamma }\cdot \delta ^{-1})\right) \text{ (}\gamma \in
\Gamma \text{)}
\end{equation*}
that is, 
\begin{equation*}
(x,\gamma )\in M\cdot N\Longleftrightarrow \exists \gamma _{1},\gamma
_{2}\in \Gamma \text{ with }(x,\gamma _{1})\in M\text{, }(x\gamma
_{1},\gamma _{2})\in N\text{ and }\gamma _{1}\gamma _{2}=\gamma \text{ }
\end{equation*}

Let $M^{1}=\overline{M}$ and for $k>1$ let 
\begin{equation*}
M^{k}=M^{k-1}\cup (M^{k-1}\cdot \overline{M})
\end{equation*}
that is, $(x,\gamma )\in M^{k}$, if and only if there exists $l\leq k$ and $%
\gamma _{1},\gamma _{2},\ldots ,\gamma _{l}\in \Gamma $ such that $\gamma
_{1}\gamma _{2}\cdots \gamma _{l}=\gamma $ and 
\begin{equation*}
(x\gamma _{1}\gamma _{2}\cdots \gamma _{i},\gamma _{i+1})\in M\text{ or }%
(x\gamma _{1}\gamma _{2}\cdots \gamma _{i}\gamma _{i+1},\gamma
_{i+1}^{-1})\in M\text{\ (}i<l\text{)}
\end{equation*}
A graphing $M$ is an $L$\emph{-graphing} if 
\begin{equation*}
\bigcup_{k}M^{k}=X\times \Gamma
\end{equation*}
The \emph{product cost} of $(\Gamma ,X)$ is defined as 
\begin{equation*}
\mathrm{pcost}(\Gamma ,X)=\inf e(M)
\end{equation*}
where the infimum is taken over all $L$-graphings $M$.

A graphing $M$ is \emph{finitely supported} if $M_{\gamma }$ is empty for
all but finitely many $\gamma \in \Gamma $. The \emph{distance} of two
graphings $M,N$ is defined as 
\begin{equation*}
d(M,N)=e(M\triangle N)=\sum_{\gamma \in \Gamma }\mu (M_{\gamma }\triangle
N_{\gamma })
\end{equation*}
where $A\triangle B$ denotes the symmetric difference of $A$ and $B$.

Let us fix a base $O$ of the topology of $X$. We call a subset of $X$ \emph{%
cylindric} with respect to $O$ if it is a finite union of elements of $O$. A
graphing $M$ is cylindric if for all $\gamma \in \Gamma $ the set $M_{\gamma
}$ is cylindric.

The following lemma says that under some assumptions on $X$ and $\Gamma $,
every $L$-graphing can be approximated with open finitely supported
cylindric $L$-graphings. Note that we do not assume that $X$ is a standard
Borel space: it can also be a finite set.

\begin{lemma}
\label{aprox}Assume that $X$ is compact, $\Gamma $ is finitely generated and
it acts by homeomorphisms on $X$. Let $M$ be an $L$-graphing of finite
measure and let $O$ be a base of the topology on $X$. Then for all $%
\varepsilon >0$ there exists an open finitely supported $L$-graphing $N$
that is cylindric with respect to $O$ and such that $d(M,N)<\varepsilon $.
\end{lemma}

\noindent \textbf{Proof.} Fix a generating set $g_{1},g_{2},\dots ,g_{d}$ of 
$\Gamma $. Let $B$ denote the graphing defined by 
\begin{equation*}
B_{\gamma }=\left\{ 
\begin{array}{cc}
X & \gamma =g_{i}\text{ for some }i\leq d \\ 
\emptyset & \text{otherwise}
\end{array}
\right.
\end{equation*}
\newline
Clearly, $B$ is an $L$-graphing of measure $d$.

List the elements of $\Gamma $ as $\gamma _{1},\gamma _{2},\dots $ and the
elements of the base $O$ as $O_{1},O_{2},\dots $.

Let $K\supseteq M$ be an open graphing satisfying $e(K\backslash
M)<\varepsilon /3$. For $n\geq 0$ let the finitely supported graphing $K(n)$
be defined by 
\begin{equation*}
K(n)_{\gamma }=\left\{ 
\begin{array}{cc}
K_{\gamma } & \gamma =\gamma _{i}\text{ for some }i\leq n \\ 
\emptyset & \text{otherwise}
\end{array}
\right.
\end{equation*}
Then $K$ has finite edge-measure, implying 
\begin{equation*}
\lim_{n\rightarrow \infty }d(K(n),K)=0
\end{equation*}
For $d\geq 2$, we have 
\begin{equation*}
K^{d}=\bigcup_{n}K(n)^{d}
\end{equation*}
which yields 
\begin{equation*}
\bigcup_{n}K(n)^{n}=\bigcup_{n}K^{n}=X\times \Gamma \supseteq B.
\end{equation*}
\ 

It is easy to see that $K(n)^{n}$ is open and $K(n)^{n}\subseteq
K(n+1)^{n+1} $ ($n\geq 2$). Since $B$ is compact, there exists $k$ such
that, setting $L=K(k)$, we have $d(L,K)<\varepsilon /3$ and $L^{k}\supseteq
B $. This implies that 
\begin{equation*}
\bigcup_{n}K^{n}\supseteq \bigcup_{n}B^{n}=X\times \Gamma
\end{equation*}
that is, $L$ is an $L$-graphing.

For $n\geq 1$ let the graphing $L(n)$ be defined by 
\begin{equation*}
L(n)_{\gamma }=\bigcup_{\substack{ 1\leq i\leq n  \\ O_{i}\subseteq
L_{\gamma }}}O_{i}
\end{equation*}
Then $L(n)$ is a finitely supported cylindric graphing ($n\geq 1$) and $%
\bigcup_{n}L(n)=L$. So, using the same argument as above we have 
\begin{equation*}
\bigcup_{n}L(n)^{n}=\bigcup_{n}L^{n}=X\times \Gamma \supseteq B
\end{equation*}
and using the compactness of $B$ again, there exists $k$ such that, setting $%
N=L(k)$, we have both $N^{k}\supseteq B$ and $d(N,L)<\varepsilon /3$. So $N$
is an $L$-graphing of distance at most $\varepsilon $ from $M$ and the lemma
is proved. $\square $

\bigskip

Now we will analyze the connection between Borel subgraphs and graphings.
Every graphing $M$ defines a Borel subgraph $\Phi (M)$ of $E$ as follows: 
\begin{equation*}
\Phi (M)=\left\{ (x,x\cdot \gamma )\mid (x,\gamma )\in M\right\}
\end{equation*}

As we will see, the map $\Phi $ is surjective. Note that it is bijective if
and only if the action of $\Gamma $ on $X$ is free.

\begin{lemma}
\label{gcostcost}We have 
\begin{equation*}
\mathrm{pcost}(\Gamma ,X)\geq \mathrm{cost}(\Gamma ,X)\text{.}
\end{equation*}
If $\Gamma $ acts essentially freely on $X$, then we have equality.
\end{lemma}

\noindent \textbf{Proof.} Let $M$ be a graphing. It is easy to see that 
\begin{equation*}
\Phi (M^{k})=\Phi (M)^{k}\text{ (}k\geq 2\text{).}
\end{equation*}
This implies that if $M$ is an $L$-graphing then $\Phi (M)$ is a Borel
subgraph of $E$ that spans $E$.

Now for $x\in X$ the degree 
\begin{eqnarray}
\deg _{S}x &=&\left| \left\{ y\in X\mid \exists \gamma \in \Gamma \text{
with }x\in M_{\gamma }\text{ and }y=x\cdot \gamma \right\} \right| \leq 
\notag \\
&\leq &\left| \left\{ \gamma \in \Gamma \mid x\in M_{\gamma }\right\}
\right| =\deg _{M}(x)
\end{eqnarray}
which implies 
\begin{equation*}
e(\Phi (M))\leq e(M)
\end{equation*}
It follows that 
\begin{equation*}
\mathrm{pcost}(\Gamma ,X)\geq \mathrm{cost}(\Gamma ,X)
\end{equation*}

Assume that $\Gamma $ acts essentially freely on $X$. Let us list the
elements of $\Gamma $ as $\gamma _{1},\gamma _{2},\dots $

Let $S$ be a Borel subgraph of $E$ spanning $E$. For each $(x,y)\in S$ let
us define $f(x,y)$ to be the first element of $\Gamma $ such that 
\begin{equation*}
y=x\cdot f(x,y)
\end{equation*}
Let us define the subset $M$ by 
\begin{equation*}
M=\left\{ (x,\gamma )\in X\times \Gamma \mid f(x,x\cdot \gamma )=\gamma
\right\}
\end{equation*}
Then $M$ is a graphing that satisfies $\Phi (M)=S$. Also, for almost all $%
x\in X$ there is equality in 1) which yields 
\begin{equation*}
e(M)=e(S)
\end{equation*}
However, the graphing $M$ may not be an $L$-graphing. Let 
\begin{equation*}
N=\bigcup_{n}M^{n}
\end{equation*}
Then 
\begin{equation*}
\Phi (N)=\bigcup_{n}\Phi (M^{n})=\bigcup_{n}\Phi (M)^{n}=\bigcup_{n}S^{n}=E
\end{equation*}
Let $(x,\gamma )\in X\times \Gamma \backslash N$. Then $(x,x\cdot \gamma
)\in S$ so for $\delta =f(x,x\cdot \gamma )\neq \gamma $ we have $x\cdot
\gamma =x\cdot \delta $, implying $\gamma \delta ^{-1}\in Stab_{\Gamma }(x)$%
. Since the action of $\Gamma $ is essentially free, we obtain $e(X\times
\Gamma \backslash N)=0$. But then 
\begin{equation*}
M^{\prime }=M\cup (X\times \Gamma \backslash N)
\end{equation*}
is an $L$-graphing of measure $e(M^{\prime })=e(S)$. This implies 
\begin{equation*}
\mathrm{pcost}(\Gamma ,X)\leq \mathrm{cost}(\Gamma ,X)
\end{equation*}
so equality holds as claimed. $\square $

\section{Boundary action and rank gradient \label{prodcost}}

In this section we first introduce coset trees and boundary representations.
Then we describe the product cost of a boundary representation in terms of
the rank gradient of the chain. This allows us to prove Theorem \ref
{coststable}.

\bigskip

Let $(\Gamma _{n})$ be a chain in $\Gamma $. Then the \emph{coset tree} $%
T=T(\Gamma ,(\Gamma _{n}))$ of $\Gamma $ with respect to $(\Gamma _{n})$ is
defined as follows. The vertex set of $T$ equals 
\begin{equation*}
T=\left\{ \Gamma _{n}g\mid n\geq 0,g\in \Gamma \right\}
\end{equation*}
and the edge set is defined by inclusion, that is, 
\begin{equation*}
(\Gamma _{n}g,\Gamma _{m}h)\text{ is an edge in }T\text{ if }m=n+1\text{ and 
}\Gamma _{n}g\supseteq \Gamma _{m}h
\end{equation*}
Then $T$ is a tree rooted at $\Gamma $ and every vertex of level $n$ has the
same number of children, equal to the index $\left| \Gamma _{n}:\Gamma
_{n+1}\right| $. The right actions of $\Gamma $ on the coset spaces $\Gamma
/\Gamma _{n}$ respect the tree structure and so $\Gamma $ acts on $T$ by
automorphisms. This action is called the \emph{tree representation} of $%
\Gamma $ with respect to $(\Gamma _{n})$.

The boundary $\partial T$ of $T$ is defined as the set of infinite rays
starting from the root. The boundary is naturally endowed with the product
topology and product measure coming from the tree. More precisely, for $%
t=\Gamma _{n}g\in T$ let us define $\mathrm{Sh}(t)\subseteq \partial T$, the 
\emph{shadow} of $t$ as 
\begin{equation*}
\mathrm{Sh}(t)=\left\{ x\in \partial T\mid t\in x\right\}
\end{equation*}
the set of rays going through $t$. Set the base of topology on $\partial T$
to be the set of shadows and set the measure of a shadow to be 
\begin{equation*}
\mu (\mathrm{Sh}(t))=1/\left| \Gamma :\Gamma _{n}\right| .
\end{equation*}
This turns $\partial T$ into a totally disconnected compact space with a
Borel probability measure $\mu $. The group $\Gamma $ acts on $\partial T$
by measure-preserving homeomorphisms; we call this action the \emph{boundary
representation of }$\Gamma $ with respect to $(\Gamma _{n})$.

\begin{lemma}
The action of $\Gamma $ on $\partial T$ is ergodic and minimal (that is,
every orbit is dense).
\end{lemma}

\noindent \textbf{Proof.} Let $A\subseteq \partial T$ be a measurable $%
\Gamma $-invariant subset such that $\mu (A)>0$. Then using the Lebesgue
density theorem, for all $\varepsilon >0$ there exists $t\in T$ of level $n$
with 
\begin{equation*}
\mu (\mathrm{Sh}(t)\cap A)\geq (1-\varepsilon )\mu (\mathrm{Sh}(t))
\end{equation*}
Since $\Gamma $ acts transitively on the $n$-th level of $T$, invariance
implies the same inequaility for all $u\in T$ of level $n$. Adding up, we
get 
\begin{equation*}
\mu (A)=\sum_{u\in T\text{ of level }n}\mu (\mathrm{Sh}(u)\cap A)\geq
1-\varepsilon 
\end{equation*}
which implies $\mu (A)=1$.

Now let $x\in \partial T$ and let $t\in T$. Let $t^{\prime }\in T$ be the
vertex of the same level as $t$ contained in $x$ and let $g\in \Gamma $ with 
$t^{\prime }g=t$. Then $xg\in \mathrm{Sh}(t)$. We proved that the orbit of $%
x $ is dense in $\partial T$. $\square $

\bigskip

There are various levels of faithfulness of a boundary representation. Let 
\begin{equation*}
\partial T_{free}=\left\{ x\in \partial T\mid Stab_{\Gamma }(x)=1\right\} 
\text{.}
\end{equation*}
We say that the action is \emph{essentially free} (or that the chain is
Farber), if $\mu (\partial T\backslash \partial T_{free})=0$. The action is 
\emph{topologically free} if $\partial T\backslash \partial T_{free}$ is
meagre, i.e., a countable union of nowhere dense closed sets. The action is 
\emph{free}, if $\partial T_{free}=\partial T$. Note that the Farber
condition has been introduced by Farber in \cite{farber} in another
equivalent formulation (see also \cite{berggab} for a relevant result).

It is easy to see that the following implications hold for the action of $%
\Gamma $ on $\partial T$: 
\begin{eqnarray*}
(\Gamma _{n})\text{ is normal and }\left( \cap \Gamma _{n}=1\right)
&\Longrightarrow &\text{free}\Longrightarrow \text{essentially free}%
\Longrightarrow \\
&\Longrightarrow &\text{topologically free}\Longleftrightarrow X_{free}\neq
\emptyset \Longrightarrow \text{faithful}
\end{eqnarray*}

For all but the third arrow it is easy to find examples showing that the
reverse implications do not hold. We shall discuss these classes more in
Section \ref{examples}.

Note that a deep result of Stuck and Zimmer \cite{stuck} tells us that every
faithful ergodic measure preserving action of a higher rank semisimple real
lattice on a probability space is essentially free. In particular, every
faithful boundary representation of such a lattice is essentially free.

\bigskip

Let $\Gamma $ be a group and let $X$ be a set. A directed $\Gamma $\emph{%
-labeled} graph is a triple $(V,E,f)$ where $(V,E)$ is a directed graph with
vertex set $V$ and edge list $E$ and $f$ is a function from $E$ to $\Gamma $%
. Note that we allow multiple edges and we make no restriction on the
labeling $f$. For a directed $\Gamma $-labeled graph 
\begin{equation*}
G=(V,E=(e_{1},\ldots ,e_{n}),f)
\end{equation*}
let $U(G)=(V,U(E))$ denote the undirected graph with vertex set $V$ and edge
list 
\begin{equation*}
U(E)=(\overline{e_{1}},\ldots ,\overline{e_{n}})
\end{equation*}
where $\overline{e}$ denotes the unordered pair obtained from the ordered
pair $e$.

Let $G$ be a directed $\Gamma $-labeled graph and let $v\in V$. Then we can
define a natural map 
\begin{equation*}
\Phi _{v}:\pi _{1}(U(G),v)\rightarrow \Gamma
\end{equation*}
from the fundamental group of $U(G)$ based at $v$ to $\Gamma $ as follows.
For a loop $l=(e_{1},\ldots ,e_{k})$ in $U(G)$ starting at $v$ let 
\begin{equation*}
\Phi _{v}(l)=\prod_{i=0}^{n-1}f^{\pm 1}(e_{i})
\end{equation*}
where the sign depends on whether we travel along $e_{i}$ preserving its
original orientation in $G$ or not. The following lemma is straightforward.

\begin{lemma}
\label{homom}The map $\Phi _{v}$ is a group homomorphism.
\end{lemma}

Now we express the rank gradient of a chain in terms of the product cost of
the action on the boundary of the coset tree. Note that we do not make any
assumptions on the boundary representation. In fact, we do not even assume
that the chain is infinite!

\begin{theorem}
\label{rankgcost}Let $(\Gamma _{n})$ be any chain in $\Gamma $. Then 
\begin{equation*}
\mathrm{RG}(\Gamma ,(\Gamma _{n}))=\mathrm{pcost}(E)-1
\end{equation*}
where $E=E(\partial T(\Gamma ,(\Gamma _{n})))$ denotes the orbit relation on 
$\partial T(\Gamma ,(\Gamma _{n}))$ defined by the action of $\Gamma $.
\end{theorem}

\noindent \textbf{Proof.} Let $s=\mathrm{RG}(\Gamma ,(\Gamma _{n}))$ and let 
$c=\mathrm{pcost}(E)$.

First we show $s+1\geq c$. Let $\varepsilon >0$. Then there exists $n$ such
that 
\begin{equation*}
\frac{d(\Gamma _{n})-1}{\left| \Gamma :\Gamma _{n}\right| }<s+\varepsilon
\end{equation*}
that is, $\Gamma _{n}$ can be generated by at most 
\begin{equation*}
d=\left\lfloor (s+\varepsilon )\left| \Gamma :\Gamma _{n}\right|
\right\rfloor +1
\end{equation*}
elements where $\left\lfloor x\right\rfloor $ is the floor of $x$. Let $%
h_{1},h_{2},\ldots ,h_{d}$ be such a generating set and let 
\begin{equation*}
\gamma _{1},\gamma _{2},\ldots ,\gamma _{\left| \Gamma :\Gamma _{n}\right| }
\end{equation*}
be a coset representative system for $\Gamma _{n}$ in $\Gamma $. We can
assume that $\gamma _{1}=1$. Let 
\begin{equation*}
Y=\left\{ h_{i}\mid i\leq d\right\} \cup \left\{ h_{i}^{-1}\mid i\leq
d\right\} \cup \left\{ 1\right\}
\end{equation*}
Let us define the graphing $M$ as follows: 
\begin{equation*}
M_{\gamma }=\left\{ 
\begin{array}{cc}
\mathrm{Sh}(\Gamma _{n}) & \text{if }\gamma =h_{i}\text{ (}i\geq 1\text{) or 
}\gamma =\gamma _{i}\text{ (}i>1\text{)} \\ 
\emptyset & \text{otherwise}
\end{array}
\right.
\end{equation*}
We claim that $M$ is an $L$-graphing. First, we have 
\begin{equation*}
\overline{M}_{\gamma }\supseteq \left\{ 
\begin{array}{cc}
\mathrm{Sh}(\Gamma _{n}) & \text{if }\gamma \in Y \\ 
\mathrm{Sh}(\Gamma _{n}) & \text{if }\gamma =\gamma _{i}\text{ (}i\leq
\left| \Gamma :\Gamma _{n}\right| \text{)} \\ 
\mathrm{Sh}(\Gamma _{n}\gamma _{i}) & \text{if }\gamma =\gamma _{i}^{-1}%
\text{ (}i\leq \left| \Gamma :\Gamma _{n}\right| \text{)}
\end{array}
\right.
\end{equation*}
Let $(x,\gamma )\in \partial T\times \Gamma $. Then there exists $a,b\leq
\left| \Gamma :\Gamma _{n}\right| $ such that $x\in \mathrm{Sh}(\Gamma
_{n}\gamma _{a})$ and $\gamma _{a}\gamma \gamma _{b}^{-1}\in \Gamma _{n}$.
This implies that there are elements $y_{1},y_{2},\ldots ,y_{k}\in Y$ such
that $\gamma _{a}\gamma \gamma _{b}^{-1}=y_{1}y_{2}\cdots y_{k}$. Using that 
$Y\subseteq \Gamma _{n}$ we get 
\begin{equation*}
(x,\gamma )=(x,\gamma _{a}^{-1}y_{1}y_{2}\cdots y_{k}\gamma _{b})\in 
\overline{M}^{k+2}.
\end{equation*}
Now the edge measure of $M$ equals 
\begin{equation*}
e(M)=\frac{1}{\left| \Gamma :\Gamma _{n}\right| }(d+\left| \Gamma :\Gamma
_{n}\right| -1)\leq 1+s+\varepsilon
\end{equation*}
which implies $c\leq 1+s$ as claimed.

Now we show that $s\leq c-1$ holds. Let $\varepsilon >0$. Then there exists
an $L$-graphing $M$ of $E$ of edge-measure at most $c+\varepsilon /2$. Using
Lemma \ref{aprox} (setting the set of shadows as base of the topology) there
exists a finitely supported cylindric $L$-graphing $N$ of $E$ with $%
d(M,N)<\varepsilon /2$. This implies that the edge-measure $%
e(N)<c+\varepsilon $.

Let $n$ be a natural number such that for every $\gamma \in \Gamma $ the set 
$N_{\gamma }$ is a union of shadows of some cosets of $\Gamma _{n}$. Let 
\begin{equation*}
V=\left\{ \Gamma _{n}\gamma \mid \gamma \in \Gamma \right\} \text{ and }%
v=\Gamma _{n}
\end{equation*}
Let us define the undirected $\Gamma $-labeled graph $G=(V,E,f)$ as follows.
For each $\gamma \in \Gamma $ and $w\in V$ where $\mathrm{Sh}(w)\subseteq
N_{\gamma }$ let us add the edge $(w,w\gamma )$ to the list $E$ with label $%
f(w,w\gamma )=\gamma $.

Let us consider the map $\Phi _{v}:\pi _{1}(U(G),v)\rightarrow \Gamma $. We
claim that the image of $\Phi _{v}$ is 
\begin{equation*}
\Phi _{v}(\pi _{1}(U(G),v))=\Gamma _{n}\text{.}
\end{equation*}
First, let $l=(e_{1},\ldots ,e_{k})$ be a loop in $U(G)$ starting at $v$.
Then $v\Phi _{v}(l)=v$ so $\Phi _{v}(l)\in \Gamma _{n}$.

Second, let $h\in \Gamma _{n}$. Let $x\in \mathrm{Sh}(v)$ be an arbitrary
element. Then since $N$ is an $L$-graphing, we have $(x,h)\in N^{k}$ for
some $k$. Thus there exist $\gamma _{1},\gamma _{2},\ldots ,\gamma _{k}\in
\Gamma $ such that $\gamma _{1}\gamma _{2}\cdots \gamma _{k}=h$ and 
\begin{equation*}
(x\gamma _{1}\gamma _{2}\cdots \gamma _{i},\gamma _{i+1})\in N\text{ or }%
(x\gamma _{1}\gamma _{2}\cdots \gamma _{i}\gamma _{i+1},\gamma
_{i+1}^{-1})\in N\text{\ (}i<k\text{)}
\end{equation*}
Let $x_{i}=x\gamma _{1}\gamma _{2}\cdots \gamma _{i}$ ($0\leq i\leq k$).
Then $vh=v$ so $x_{0}=v=x_{k}$. Also, for all $0\leq i<k$ there is an edge
in $G$ either from $x_{i}$ to $x_{i+1}$ labeled by $\gamma _{i+1}$ or from $%
x_{i+1}$ to $x_{i}$ labeled by $\gamma _{i+1}^{-1}$. Thus these edges form a
loop $l$ in $U(G)$ with $\Phi _{v}(l)=h$. The claim follows.

The number of vertices of $G$ equals $\left| \Gamma :\Gamma _{n}\right| $
while the number of edges of $G$ equals $e(N)\left| \Gamma :\Gamma
_{n}\right| $. Hence, the same holds for $U(G)$. Using Lemma \ref{homom} and
the formula for the rank of the fundamental group of a graph we get 
\begin{equation*}
d(\Gamma _{n})\leq d(\pi _{1}(U(G),v))=e(N)\left| \Gamma :\Gamma _{n}\right|
-\left| \Gamma :\Gamma _{n}\right| +1
\end{equation*}
which yields 
\begin{equation*}
\frac{d(\Gamma _{n})-1}{\left| \Gamma :\Gamma _{n}\right| }\leq
e(N)-1<c-1+\varepsilon \text{.}
\end{equation*}
This shows that $s\leq c-1$.

The theorem is proved. $\square $

\bigskip

Now Theorem \ref{coststable} follows immediately.

\bigskip

\noindent \textbf{Proof of Theorem \ref{coststable}.} Using Theorem \ref
{rankgcost} and Lemma \ref{gcostcost} we get 
\begin{equation*}
\mathrm{RG}(\Gamma ,(\Gamma _{n}))=\mathrm{pcost}(E)-1=\mathrm{cost}(E)-1
\end{equation*}
as claimed. $\square $

\section{Applications and examples \label{examples}}

In this section we introduce the absolute rank gradient of a group and
compute it for some important classes of groups. We will later use these
results in Section \ref{manifolds}. Then we discuss what happens if we relax
the Farber condition.

Let us define the \emph{absolute rank gradient} of $\Gamma $ as 
\begin{equation*}
\mathrm{RG}(\Gamma )=\inf_{H}\frac{d(H)-1}{\left| \Gamma :H\right| }
\end{equation*}
where $H$ runs through all subgroups of $\Gamma $ of finite index.

Let $H,K\leq \Gamma $ be subgroups of finite index with $H\leq K$. Using the
Nielsen-Schreier theorem on $H$ and $K$, we get $d(H)-1\leq \left|
K:H\right| (d(K)-1)$ which yields 
\begin{equation*}
\frac{d(H)-1}{\left| \Gamma :H\right| }\leq \frac{d(K)-1}{\left| \Gamma
:K\right| }
\end{equation*}
So for any chain $(\Gamma _{n})$ in $\Gamma $, the sequence $(d(\Gamma
_{n})-1)/\left| \Gamma :\Gamma _{n}\right| $ is non-increasing and the
definition of $\mathrm{RG}(\Gamma ,(\Gamma _{n}))$ makes sense. In fact, as
the authors show in \cite{abnik}, if the sequence stabilizes then $\Gamma $
is virtually free.

We shall make use of the following easy lemma.

\begin{lemma}
\label{belerak}Let $\Gamma $ be a finitely generated, residually finite
group and let $H\leq \Gamma $ be a subgroup of finite index. Then there
exists a chain $(\Gamma _{n})$ in $\Gamma $ such that $\Gamma _{1}=H$ and
the boundary representation of $\Gamma $ with respect to $(\Gamma _{n})$ is
free.
\end{lemma}

\noindent \textbf{Proof.} Let $K$ be the core of $H$, that is, $%
K=\bigcap_{g\in \Gamma }H^{g}$. Then $K$ is a normal subgroup of $\Gamma $
of finite index. For $n\geq 1$ let 
\begin{equation*}
\Delta _{n}=\bigcap_{L\leq \Gamma \text{ of index }n}L
\end{equation*}
Since $\Gamma $ has only finitely many subgroups of a given index, each $%
\Delta _{n}$ is a normal subgroup of finite index in $\Gamma $. Also, since $%
\Gamma $ is residually finite, $\cap _{n}\Delta _{n}=1$. Now let $\Gamma
_{0}=\Gamma $, let $\Gamma _{1}=H$ and for $n\geq 2$ let $\Gamma _{n}=K\cap
\Delta _{n}$. Then $\cap _{n}\Gamma _{n}=1$ and the chain $(\Gamma _{n})$
consists of normal subgroups of $\Gamma $ (except for $n=1$). Let $g\in
\Gamma $ with $g\neq 1$. Then there exists $n>1$ such that $g\notin \Gamma
_{n}$. This implies that $g$ acts fixed point freely on $\Gamma /\Gamma _{n}$
and thus also on the boundary $\partial T(\Gamma ,(\Gamma _{n}))$. In other
words, the action of $\Gamma $ on $\partial T$ is free. $\square $

\bigskip

Theorem \ref{coststable} now gives us the following on the absolute rank
gradient.

\begin{corollary}
\label{absolute}Let $\Gamma $ be a finitely generated, residually finite
group. Then 
\begin{equation*}
\mathrm{RG}(\Gamma )=\mathrm{cost}(\Gamma ,\widehat{\Gamma })-1
\end{equation*}
where $\widehat{\Gamma }$ denotes the profinite completion of $\Gamma $.
\end{corollary}

\noindent \textbf{Proof.} Let $\Delta _{n}$ be as in the proof of Lemma \ref
{belerak}. Then for all $H\leq \Gamma $ of index $n$, we have $\Delta
_{n}\leq H$. Hence $\widehat{\Gamma }_{(\Delta _{n})}$, the profinite
completion of $\Gamma $ with respect to $(\Delta _{n})$ equals $\widehat{%
\Gamma }$ and by Theorem \ref{coststable} we get 
\begin{equation*}
\mathrm{RG}(\Gamma )=\mathrm{RG}(\Gamma ,\Delta _{n})=\mathrm{cost}(\Gamma ,%
\widehat{\Gamma }_{(\Delta _{n})})-1=\mathrm{cost}(\Gamma ,\widehat{\Gamma }%
)-1
\end{equation*}
as claimed. $\square $

\bigskip

\noindent \textbf{Proof of Corollary \ref{amennormal}.} Let $\Gamma $ be a
finitely generated residually finite group with an infinite amenable normal
subgroup. Let $(\Gamma _{n})$ be a Farber chain in $\Gamma $ -- such chain
exists by Lemma \ref{belerak}. Let $E=E(\partial T(\Gamma ,(\Gamma _{n})))$.
Then by \cite{gabor} $\Gamma $ has fixed price $1$ and so, using Theorem \ref
{coststable} we have 
\begin{equation*}
\mathrm{RG}(\Gamma )=\mathrm{RG}(\Gamma ,(\Gamma _{n}))=\mathrm{cost}(E)-1=0.
\end{equation*}
$\square $

\bigskip

Note that in \cite{abnik} we present an interesting alternative
combinatorial proof of the corollary in the case when $\Gamma $ itself is
infinite amenable. However, we have no proof for the general case that does
not use analysis.

\bigskip

\noindent \textbf{Ascending HNN\ extensions.} Let $A$ be a finitely
generated group and $f:A\rightarrow A$ a homomorphism. Then the ascending
HNN extension with base $A$ is the group $\Gamma $ with presentation 
\begin{equation*}
\left\langle A,t\ |\ a^{t}=f(a)\ \forall a\in A\right\rangle
\end{equation*}
For example every free by cyclic group has such a presentation.

For each $n\in \mathbb{N}$ let $\Gamma _{n}=\langle A,t^{n}\rangle \leq
\Gamma $. It is easy to see that $\Gamma _{n}$ is a normal subgroup of index 
$n$ in $\Gamma $. As $d(\Gamma _{n})\leq d(A)+1$ we get the following
result, first observed by Lackenby in \cite{lack}.

\begin{proposition}
\label{hnn} If the group $\Gamma $ is a residually finite ascending HNN
extension, then $\mathrm{RG}(\Gamma )=0$.
\end{proposition}

For instance, free by cyclic groups are residually finite and so they have
rank gradient $0$. In fact, more generally we have the following.

\begin{proposition}
\label{fgnormal} Let $\Gamma $ be finitely generated residually finite group
which has a finitely generated normal subgroup $N$ such that $\Gamma /N$ has
subgroups of arbitrarily large index (i.e. the profinite completion $%
\widehat{\Gamma /N}$ of $\Gamma /N$ is infinite). Then $\mathrm{RG}(\Gamma
)=0$.
\end{proposition}

\noindent \textbf{Proof.} Using the assumptions one can find a sequence of
subgroups $H_{i}<\Gamma $ such that both 
\begin{equation*}
a_{i}=[\Gamma :NH_{i}]\text{ and }b_{i}=[N:N\cap H_{i}]
\end{equation*}
tend to infinity as $i$ tends to infinity. Now $[\Gamma :H_{i}]=a_{i}b_{i}$, 
$d(N\cap H_{i})\leq b_{i}d(N)$ and $d(NH_{i}/N)\leq a_{i}d(\Gamma )$. Hence
from $d(H_{i})\leq d(N\cap H_{i})+d((NH_{i}/N)$ we obtain 
\begin{equation*}
\frac{d(H_{i})-1}{[\Gamma :H_{i}]}<\frac{d(\Gamma )}{b_{i}}+\frac{d(N)}{a_{i}%
}
\end{equation*}
which tends to $0$ as $i$ tends to infinity. This implies $\mathrm{RG}%
(\Gamma )=0$. $\square $

\bigskip

The result that the groups in Proposition \ref{fgnormal} have a measurable
action with cost $1$, has been proved by Gaboriau \cite{gabor}, even without
assuming residual finiteness of $\Gamma $ and $|\widehat{\Gamma /N}|=\infty $%
, but that does not give anything on the rank gradient. We note that the
condition that $\widehat{\Gamma /N}$ is infinite seems quite mild. 

\bigskip

Now we discuss rank gradient of chains that are not Farber. As we shall see,
the situation can be quite different.

Our first example is the so-called \emph{lamplighter group}, the wreath
product 
\begin{equation*}
\Gamma =\mathbb{Z}/2\mathbb{Z}\ \mathrm{wr}\ \mathbb{Z}\text{.}
\end{equation*}
This group is metabelian, in particular, it is amenable, so Corollary \ref
{amennormal} implies that the rank gradient of any Farber chain is $0$.
However, consider the canonical surjections 
\begin{equation*}
\phi _{n}:\mathbb{Z}/2\mathbb{Z}\text{ }\mathrm{wr}\ \mathbb{Z}\rightarrow 
\mathbb{Z}/2\mathbb{Z}\ \mathrm{wr}\ \mathbb{Z}/2^{n}\mathbb{Z}\text{ and }%
\pi _{n}:\mathbb{Z}/2\mathbb{Z}\text{ }\mathrm{wr}\text{ }\mathbb{Z}/2^{n}%
\mathbb{Z}\rightarrow \mathbb{Z}/2^{n}\mathbb{Z}
\end{equation*}
and let $\Gamma _{n}=\mathrm{Ker}(\phi _{n}\pi _{n})$. Then $(\Gamma _{n})$
is a normal chain, $\Gamma /\Gamma _{n}\approxeq \mathbb{Z}/2^{n}\mathbb{Z}$
and $\phi _{n}(\Gamma _{n})\approxeq (\mathbb{Z}/2\mathbb{Z})^{2^{n}}$ which
implies $d(\Gamma _{n})\geq 2^{n}$. This yields 
\begin{equation*}
\mathrm{RG}(\Gamma ,(\Gamma _{n}))\geq 1.
\end{equation*}

In the above example the chain is normal but with nontrivial intersection,
and so the boundary representation is not faithful. We do not know whether
there exists an amenable group $\Gamma $ and a chain $(\Gamma _{n})$ such
that the boundary representation is faithful, but the chain has positive
rank gradient.

Our last example is a slight variation of one discussed by Bergeron and
Gaboriau in \cite{berggab}. The difference is that they estimate the first
Betti number where we estimate the rank.

\begin{proposition}
\label{nn}There exist a virtually free group $\Gamma $ and an interval $%
[x,y)\subset \mathbb{R}$ such that for every $\alpha \in \lbrack x,y)$ there
exists a subnormal chains of subgroups $\Gamma =H_{0}>H_{1}>H_{2}>\cdots $
with trivial intersection, such that $\mathrm{RG}(\Gamma ,(H_{i}))=\alpha $.

Moreover in this situtation we can have the strict inequalities 
\begin{equation*}
\mathrm{RG}(\Gamma ,(H_{i}))>\lim_{i\rightarrow \infty }\frac{b_{1,p}(H_{i})%
}{\left| \Gamma :H_{i}\right| }>\lim_{i\rightarrow \infty }\frac{b_{1}(H_{i})%
}{\left| \Gamma :H_{i}\right| },
\end{equation*}
where for a prime $p$ the integer $b_{1,p}(H)=$ $\dim _{\mathbb{F}%
_{p}}H/[H,H]H^{p}$ is the $p$-homology of $H$.
\end{proposition}

In particular, these chains are not Farber.

\bigskip

\noindent \textbf{Proof of Proposition \ref{nn}. }Let $\Gamma =A\ast \mathbb{%
Z}$, where $A$ is a finite group. Viewing $\Gamma $ as a trivial HNN
extension we see that $\Gamma $ has a right transitive action on both the
vertices and the edges of a regular $2|A|$-valent tree $T$ (so that the
quotient $I_{0}$ is a single vertex with a looped edge). We can direct the
edges of $T$ so that $\Gamma $ acts preserving this orientation. There are
exactly $a=|A|$ in-edges and out-edges from every vertex of $T$. The action
of $\Gamma $ is regular on the edges of $T$. We fix an edge $e_{0}$ of $T$
and a vertex $v_{0}$ at one end of $e_{0}$ so that the stabilizer of $v_{0}$
is $A.\bigskip $

Suppose that $H$ is a subgroup of $\Gamma $ of index $n$. Then the quotient $%
I=$ $T/H$ is a finite graph with a covering map $p:T\rightarrow I=T/H$ and $%
H $ can be recovered from $p$ as $\{g\in \Gamma \ |\ p(e_{0})=p(e_{0}g)\}$.

The vertices of $I$ are in 1-1 correspondence with the double cosets $%
A\backslash \Gamma /H$. Given a vertex $v\in $ $I$ define $S_{v}=$ $%
A^{g}\cap H$. where $g$ is a fixed representative of $A\backslash \Gamma /H$
such that $p(v_{0}g)=v$. There are exactly $n$ edges in $I$ and we have $%
n=\sum_{v\in I}[A:S_{v}]$.

The Bass-Serre theory determines the structure of $H$ as follows: 
\begin{equation*}
H=F\ast (\ast _{v\in I}S_{v})
\end{equation*}
where the free group $F$ is the fundamental group of $I$. The group $F$ has
rank $n-p+1$, where $p$ is the number of vertices in $I$.

Let $X$ be the set of vertices of $I$ with out-valency $1$ and let $%
Y=I\backslash X$ be the set of vertices with out-valency $a=|A|$. Let $%
\left| X\right| =\mu p$. Then we have 
\begin{equation*}
n=p(\mu +(1-\mu )|A|),\text{ hence }p=\frac{n}{\mu +(1-\mu )|A|}.
\end{equation*}
We get that exactly $\mu p$ of the vertex stabilizers $S_{v}$ are equal to $%
A $ and the rest are trivial. From the Grushko--Neumann theorem (see \cite
{gru} and \cite{neu}) it now follows that the rank of $H$ is 
\begin{equation*}
d(H)=n-p+\mu pd(A)+1=n\left( 1+\frac{\mu d(A)-1}{\mu +(1-\mu )|A|}\right) +1
\end{equation*}
while $\beta _{1}(H)=n-np+1.$

Now suppose that we have a sequence of finite oriented graphs with covering
maps 
\begin{equation*}
I_{0}\leftarrow I_{1}\leftarrow I_{2}\leftarrow \cdots \leftarrow T
\end{equation*}
with the following two properies: \medskip

1. Each vertex of each $I_{j}$ has valency $2$ or $2a$. The proportions $\mu
_{j}=|X_{j}|/|I_{j}|$ of the vertices of $I_{j}$ with out- and in-valency $1$
form a decreasing sequence which tends to some prescribed limit $\mu =\mu
_{\infty }\in \lbrack 0,1)$. \medskip

2. Given any integer $k$ let $B_{k}$ be the ball of radius $k$ in the $2a$%
-regular oriented tree $T$ centered at the vertex $v_{0}$. Then there is an
integer $m$ such that the covering map $T\rightarrow I_{m}$ is injective on $%
B_{k}.$

\bigskip

Then the groups $H_{j}$ which correspond to the graphs $I_{j}$ form a chain
in $\Gamma $ with trivial intersection and such that 
\begin{equation*}
\lim_{i\rightarrow \infty }\frac{d(H_{i})-1}{\left| \Gamma :H_{i}\right| }=1+%
\frac{\mu d(A)-1}{\mu +(1-\mu )|A|}
\end{equation*}

The first part of Proposition \ref{nn} now follows by setting $x=1-\frac{1}{%
|A|}$ (when $\mu _{\infty }=0$) and $y=d(A)$ (when $\mu _{\infty
}\rightarrow 1$).

The existence of the graphs $I_{j}$ with properties 1 and 2 above is
essentially proved by Bergeron and Gaboriau in \cite[Section 4]{berggab} in
the case of a free product of two residually finite groups.

On the other hand, we have 
\begin{equation*}
\lim_{i\rightarrow \infty }\frac{\beta _{1}(H_{i})}{\left| \Gamma
:H_{i}\right| }=1-\frac{1}{\mu +(1-\mu )|A|}
\end{equation*}
and similarly we have 
\begin{equation*}
b_{1,p}(H)=n\left( 1+\frac{\mu b_{1,p}(A)-1}{\mu +(1-\mu )|A|}\right) +1%
\text{.}
\end{equation*}
Therefore 
\begin{equation*}
\lim_{i\rightarrow \infty }\frac{\beta _{1,p}(H_{i})}{\left| \Gamma
:H_{i}\right| }=1+\frac{\mu \beta _{1,p}(A)-1}{\mu +(1-\mu )|A|}.
\end{equation*}

So, if we choose the finite group $A$ so that $d(A)>b_{1,p}(A)>0$ and choose
the chain $(H_{i})$ with limiting ratio $\mu _{\infty }>0$ this gives 
\begin{equation*}
\mathrm{RG}(\Gamma ,(H_{i}))>\lim_{i\rightarrow \infty }\frac{b_{1,p}(H_{i})%
}{\left| \Gamma :H_{i}\right| }>\lim_{i\rightarrow \infty }\frac{b_{1}(H_{i})%
}{\left| \Gamma :H_{i}\right| }
\end{equation*}
as promised. $\square $

\section{Lattices and 3-manifolds \label{manifolds}}

In this section we further discuss the rank gradient, with a special
emphasis on lattices in Lie groups. We also define the Heegaard genus and
derive the contradiction between the Rank vs. Heegaard genus conjecture and
Fixed Price.

The following lemma follows immediately from \cite[Theorem 3]{gabor} saying
that if $H\leq \Gamma $ of finite index, then 
\begin{equation*}
\mathrm{cost}(H)-1=\left| \Gamma :H\right| (\mathrm{cost}(\Gamma )-1)
\end{equation*}
where $\mathrm{cost}(\Gamma )$ is the infimum of costs of all measurable
essentially free actions of $\Gamma $.

\begin{lemma}
\label{fixedmulti}An affirmative answer to the Fixed Price problem implies
an affirmative answer to the Multiplicativity of cost-1 problem.
\end{lemma}

\noindent \textbf{Proof.} Let $\Gamma $ be a finitely generated group acting
on $(X,\mu )$ by measure preserving maps and let $H$ be a subgroup of $%
\Gamma $ of finite index. Then by our assumption we have 
\begin{equation*}
\mathrm{cost}(H,X)-1=\mathrm{cost}(H)-1=\left| \Gamma :H\right| (\mathrm{cost%
}(\Gamma )-1)=\left| \Gamma :H\right| (\mathrm{cost}(\Gamma ,X)-1)
\end{equation*}
$\square $

\bigskip

Now we proceed to the independence of the rank gradient.

\begin{theorem}
\label{independence}Assume that the Multiplicativity of cost-1 problem has
an affirmative solution. Let $(\Gamma _{n})$ be a Farber chain in $\Gamma $.
Then 
\begin{equation*}
\mathrm{RG}(\Gamma ,(\Gamma _{n}))=\mathrm{RG}(\Gamma )
\end{equation*}
In particular, any two Farber chains have the same rank gradient in $\Gamma $%
.
\end{theorem}

\noindent \textbf{Proof.} Let $\Gamma $ act on a standard Borel probability
space $X$ essentially freely. Then using the multiplicativity assumption, we
have 
\begin{equation*}
\mathrm{cost}(\Gamma ,X)-1=\frac{\mathrm{cost}(\Gamma _{n},X)-1}{\left|
\Gamma :\Gamma _{n}\right| }\leq \frac{d(\Gamma _{n})-1}{\left| \Gamma
:\Gamma _{n}\right| }
\end{equation*}
which implies 
\begin{equation*}
\mathrm{cost}(\Gamma ,X)-1\leq \mathrm{RG}(\Gamma ,(\Gamma _{n}))\text{.}
\end{equation*}

Now if $(\Lambda _{n})$ is another Farber chain in $\Gamma $, then using
Theorem \ref{coststable} and the above inequiality for both chains and
actions, we get 
\begin{eqnarray*}
\mathrm{RG}(\Gamma ,(\Gamma _{n})) &=&\mathrm{cost}(\Gamma ,E(\Gamma
,(\Gamma _{n})))-1\leq \\
&\leq &\mathrm{RG}(\Gamma ,(\Lambda _{n}))=\mathrm{cost}(\Gamma ,E(\Gamma
,(\Lambda _{n})))-1\leq \mathrm{RG}(\Gamma ,(\Gamma _{n}))
\end{eqnarray*}
so equality holds everywhere.

Let $\varepsilon >0$ and let $H\leq \Gamma $ be a subgroup of finite index
such that 
\begin{equation*}
\frac{d(H)-1}{\left| \Gamma :H\right| }<\mathrm{RG}(\Gamma )+\varepsilon
\end{equation*}
Using Lemma \ref{belerak} there exists a Farber chain $(\Lambda _{n})$ in $%
\Gamma $ with $\Lambda _{1}=H$. This implies 
\begin{equation*}
\mathrm{RG}(\Gamma ,(\Gamma _{n}))=\mathrm{RG}(\Gamma ,(\Lambda _{n}))\leq 
\frac{d(H)-1}{\left| \Gamma :H\right| }<\mathrm{RG}(\Gamma )+\varepsilon
\end{equation*}
and so $\mathrm{RG}(\Gamma ,(\Gamma _{n}))\leq \mathrm{RG}(\Gamma )$.

On the other hand, $\mathrm{RG}(\Gamma ,(\Gamma _{n}))\geq \mathrm{RG}%
(\Gamma )$ by definition, so the theorem holds. $\square $

\bigskip

Actually, the above argument shows that assuming the multiplicativity of
cost-1, a finitely generated residually finite group has fixed price $1$ if
and only if its absolute rank gradient is $0$.

This could be relevant to decide whether uniform lattices have fixed price $%
1 $. Let $G$ be a semisimple Lie group of $\mathbb{R}$-rank at least $2$ and
let $\Gamma $ be a lattice in $G$. Gaboriau \cite{gabor} shows that if $%
\Gamma $ is non-uniform, i.e., $G/\Gamma $ is not compact, then $\Gamma $
has fixed price $1$. Since, by a theorem of Borel, for every uniform higher
rank lattice, there is a non-uniform lattice in the same ambient group,
every uniform lattice in $G$ has a measurable action of cost $1$. But fixed
price is not known for these lattices.

\begin{conjecture}
All lattices in higher rank Lie groups have absolute rank gradient $0$.
\end{conjecture}

Note that from a result of Sharma and Venkataramana \cite{venky}, every
non-uniform lattice $\Gamma $ in $G$ contains a subgroup of finite index
generated by just $3$ elements. This trivially implies $\mathrm{RG}(\Gamma
)=0$.

\bigskip

Now we discuss rank $1$ lattices. First, by \cite{gabor}, every lattice $%
\Gamma $ in $\mathrm{SL}(2,\mathbb{R})$ (or $\mathrm{PSL}(2,\mathbb{R})$)
has fixed price greater than $1$. Using Theorem \ref{coststable} we get that 
$\mathrm{RG}(\Gamma )>0$. Next we consider three examples, two of which have
been studied extensively in the literature.

\bigskip

\textbf{Example A.} Let $\Gamma $ be the Picard group $\mathrm{PSL}(2,%
\mathbb{Z}[i])$. It is known that $\Gamma $ is virtually a free by cyclic
group. There is a subgroup $T$ of index $24$ in $\Gamma $ which is the
fundamental group of the complement of the Borromean links and so (cf. \cite
{menzel}) $T$ is an extension of the free group on $4$ generators $F_{4}$ by 
$\mathbb{Z}$. From Proposition \ref{fgnormal} it follows that $\mathrm{RG}%
(\Gamma )=0$.

\bigskip

\textbf{Example B. }Let $\Gamma $ be the Bianchi group $\mathrm{PSL}(2,%
\mathbb{Z}[w])$ where $w=\exp (2\pi i/3)$. Here there is a subgroup $U$ of
index $12$ of $\Gamma $ which is the fundamental group of the figure eight
knot complement. Again, it is well known (see \cite{eightknot}) that $U$ is
the free by cyclic group with a presentation 
\begin{equation*}
\langle a,b,t\ |\ a^{t}=b^{-1},\ b^{t}=b^{2}ab\rangle .
\end{equation*}
So in this case again $\Gamma $ is virtually a free by cyclic group and so $%
\mathrm{RG}(\Gamma )=0$.

\bigskip

\textbf{Example C. }Both Examples A and B are non-uniform. It is much harder
to find uniform lattices of $\mathrm{PSL}(2,\mathbb{C})$ with rank gradient
zero. One such group has been found by Reid \cite[Section 4, Theorem 1]{AR}.
His example is a uniform arithmetic lattice $\Lambda $ (that is, $\Lambda $
is commensurable with the group of norm $1$ elements of a suitable
quaternion algebra over a number field $K$ with just one complex embedding).
Reid proves that the manifold $\mathbb{H}^{3}/\Lambda $ has a finite cover
which fibres over the circle. Group theoretrically this means that $\Lambda $
has a subgroup of finite index which is an ascending HNN extension with a
finitely generated base group $A$. Proposition \ref{hnn} now gives that the
rank gradient of $\Lambda $ is $0$.

\bigskip

\noindent \textbf{On a question of Kechris and Miller. }Examples A and B
present a negative answer to a question of Kechris and Miller \cite[Problem
35.7]{kechmill}. They asked whether for a countable infinite Euclidean
domain $D$ the group $\mathrm{SL}(2,D)$ has cost $1$ if and only if $D$ has
infinitely many units. In fact, we get that these groups have rank gradient $%
0$ and hence cost $1$. We suspect that this is the general behaviour.
Indeed, a well-known conjecture by Thurston asserts that every hyperbolic $3$%
-manifold virtually fibers over the circle. This would imply that any
lattice in $\mathrm{SL}(2,\mathbb{C})$ has cost $1$.

\bigskip

Let $M$ be a finite closed orientable $3$-manifold. A Heegaard decomposition
(or splitting) for $M$ is an expression of $M$ as a union of two isomorphic
handlebodies $h_{1},h_{2}$ of genus $g\geq 0$ with boundary surfaces $%
\partial h_{i}$ identified via a homeomorphism $f:\partial h_{1}\rightarrow
\partial h_{2}$. Such decomposition exists by \cite{JS}, Section 8.3.

The Heegaard genus $g(M)$ of $M$ is the minimal genus $g$ of the surfaces $%
\partial h_{i}$ in some Heegaard decomposition for $M$. Let $r(M)=d(\pi
_{1}(M))$ be the rank of $M$.

It is easy to see that the fundamental group $\pi _{1}(h_{i})$ surjects onto 
$\pi _{1}(M)$ and so $g(M)\geq r(M)$. In \cite{W} Waldhausen asked if there
is equality. This was proved false for Seifert manifolds by Boileau and
Zieschang in \cite{BZ}. Further work has been done by Schultens and Weidman
who construct manifolds $M$ with $g(M)=4n$ and $r(M)\leq 3n$ ($n\geq 1$).
However, the question of Waldhausen remained open for hyperbolic $3$%
-manifolds. Also, until now, for all the known counterexamples the ratio $%
g(M)/r(M)$ was at most $3/2$. Theorem \ref{manifolds} addresses both these
problems by proving that the ratio can get arbitrarily large even for
arithmetic hyperbolic\textbf{\ }$3$-manifolds\textbf{. }

\bigskip

We are ready to show that the Rank vs Heegaard genus conjecture and Fixed
price problem conflict each other.

\bigskip

\noindent \textbf{Proof of Theorem \ref{contradiction}.} Assume that the
Fixed Price problem has an affirmative solution. Hence, by Lemma \ref
{fixedmulti} and Theorem \ref{independence}, any Farber chain has rank
gradient equal to the absolute rank gradient of $\Gamma $.

Let us take the group $\Lambda $ in Example C above. Take a chain $\Lambda
=N_{0}>N_{1}>N_{2}\cdots $ of \emph{congruence} normal subgroups of $\Gamma $
with trivial intersection. Define the manifolds $M_{i}=\mathbb{H}^{3}/N_{i}$%
. So we have 
\begin{equation*}
\pi _{1}(M_{i})=N_{i}.
\end{equation*}

From results of Sarnak and Xue \cite{sarnak} (see also \cite[p. 445, Example
(f)]{lub}) it follows that if $L$ is an arithmetic lattice of $\mathrm{PSL}%
(2,\mathbb{C})$, then $L$ has property $(\tau )$ with respect to its
congruence subgroups. In particular, $\Lambda $ has property $(\tau )$ with
respect to the chain $\{N_{i}\}$ and hence by a recent result of Marc
Lackenby, \cite[Theorem 1.5 and Corollary 1.6]{lackhaken} we have 
\begin{equation*}
\lim_{i}\frac{g(M_{i})}{[\Lambda :N_{i}]}>0.
\end{equation*}

On the other hand, we have shown that $\mathrm{RG}(\Gamma )=0$, so by
Theorem \ref{independence} we have 
\begin{equation*}
\mathrm{RG}(\Gamma ,(N_{i}))=\lim_{i}\frac{r(M_{i})}{[\Lambda :N_{i}]}=0%
\text{.}
\end{equation*}
This implies that the ratio $g(M_{i})/r(M_{i})$ tends to infinity with $i$,
thus proving the theorem. $\square $

\bigskip

\noindent \textbf{Remark 1.} If the manifold $M$ is noncompact then its
canonical compactification $M^{c}$ is a finite $3$-manifold with boundary.
In this case a Heegaard decomposition of $M$ is its expression as a union of
two \emph{compression bodies} $B_{1},B_{2}$ with their positive boundaries $%
\partial B_{i+}$ (both orientable surfaces of genus $g$) identified. The
boundary of $M$ is the negative boundary $\partial B_{1-}\cup \partial
B_{2-} $. See \cite{lackhaken}, Section 3 for more details. The Heegaard
genus $g(M) $ of $M$ is defined to be the minimal genus of the surfaces
indentified in some decomposition for $M^{c}$ as above. Again it is easy to
see that in this case we have the inequality $2g(M)\geq r(M)$. Assuming the
independence of the rank gradient on the chain, for Examples A or B, just as
for compact manifolds, the ratio $g(M)/r(M)$ can be arbitrarily large.

\bigskip

\noindent \textbf{Remark 2. }Note that one does not need to obtain the
multiplicativity of cost-1 for general measurable actions to obtain the
independence of rank gradient from the chain. Trivially, it would be enough
to settle this for profinite actions. Less trivially, by the recent work of
the first author and Elek \cite{abelek}, it would be enough to settle the
multiplicativity for the standard Bernoulli action $\{0,1\}^{\Gamma }$.

\end{document}